\documentclass[a4paper,11pt,reqno,intlimits]{extarticle}
\usepackage[cp1251]{inputenc}
\usepackage[T2A]{fontenc}
\usepackage[english,russian]{babel}
\usepackage{amsfonts,amsmath,amsxtra,amsthm,amssymb,latexsym,mathrsfs,graphicx,setspace}
\usepackage{multirow}
\usepackage{url}

\usepackage{algorithm}
\usepackage{algorithmic}
\makeatletter
\renewcommand{\ALG@name}{Алгоритм}
\makeatother

\theoremstyle{plain} 
\newtheorem{theorem}{Теорема}
\newtheorem{lemma}{Лемма}

\newtheorem{corollary}{Следствие}

\theoremstyle{definition} 
\newtheorem{remark}{Замечание}

\numberwithin{equation}{section} \numberwithin{theorem}{section}
\numberwithin{lemma}{section} \numberwithin{proposition}{section}
\numberwithin{corollary}{section} \numberwithin{remark}{section}
\numberwithin{definition}{section} \numberwithin{example}{section}

\begin{document}

\renewcommand{\thesection}{}
\renewcommand{\thesubsection}{}
\renewcommand{\thetheorem}{\arabic{theorem}}
\renewcommand{\theequation}{\arabic{equation}}
\renewcommand{\thedefinition}{\arabic{definition}}
\renewcommand{\thecorollary}{\arabic{corollary}}
\renewcommand{\thelemma}{\arabic{lemma}}
\renewcommand{\theexample}{\arabic{example}}
\renewcommand{\theremark}{\arabic{remark}}

\author{Стонякин Ф.С., Alkousa M., Степанов А.Н., Баринов М.А.}
\title{\textbf{
Адаптивные алгоритмы зеркального спуска в задачах выпуклого программирования с липшицевыми ограничениями
\\
Adaptive algorithms for mirror descent in convex programming problems with Lipschitz constraints
}}

\maketitle

\begin{abstract}
Работа посвящена новым модификациям недавно предложенных адаптивных методов зеркального спуска для задач выпуклой минимизации в случае нескольких выпуклых функциональных ограничений. Предложены адаптивные методы зеркального спуска для задач двух типов. Первый тип~--- задачи с липшицевым (вообще говоря, негладким) целевым функционалом. Второй тип~--- задачи с липшицевым градиентом целевого функционала. Рассматривается также случай негладкого целевого функционала, равного максимуму гладких функционалов с липшицевым градиентом. Во всех случаях функциональные ограничения считаются выпуклыми, липшицевыми и, вообще говоря, негладкими. Предлагаемые методы позволяют сэкономить время работы алгоритма за счёт рассмотрения не всех функциональных ограничений на непродуктивных шагах. Получены оценки на скорость сходимости рассматриваемых методов. Эти оценки демонстрируют оптимальность методов с точки зрения нижних оракульных оценок. Приведены результаты численных экспериментов, иллюстрирующие преимущества предлагаемой методики для некоторых примеров.
\end{abstract}

\textbf{Ключевые слова:} адаптивный метод зеркального спуска, липшицев\\ функционал, липшицев градиент, продуктивный шаг, непродуктивный шаг.

\selectlanguage{english}
\begin{abstract}
The paper is devoted to new modifications of recently proposed adaptive methods of Mirror Descent for convex minimization problems in the case of several convex functional constraints. Methods for problems of two classes are considered. The first type of problems with Lipschitz-continuous objective (generally speaking, nonsmooth) functional. The second one is for problems with a Lipschitz-continuous gradient of the objective smooth functional. We consider the class of problems with a non-smooth objective functional equal to the maximum of smooth functionals with a Lipschitz-continuous gradient. Note that functional constraints, generally speaking, are non-smooth and Lipschitz-contionuous. The proposed modifications allow saving the algorithm running time due to consideration of not all functional constraints on non-productive steps. Estimates for the rate of convergence of the methods under consideration are obtained. The methods proposed are optimal from the point of view of lower oracle estimates. The results of numerical experiments illustrating the advantages of the proposed procedure for some examples are given.
\end{abstract}
\selectlanguage{russian}

Lipschitz-continuous gradient, productive step, non-productive step.\\


\subsection*{Введение}

Задачи минимизации негладкого функционала c ограничениями возникают в широком классе проблем современной large-scale оптимизации и её приложений \cite{bib_ttd,bib_Shpirko}. Для таких задач имеется множество методов, среди которых можно отметить метод зеркального спуска (МЗС) \cite{beck2003mirror,nemirovsky1983problem} для задач выпуклой минимизации.

Отметим, что в случае негладкого целевого функционала или функциональных ограничений естественно использовать субградиентные методы, восходящие к хорошо известным работам \cite{polyak1967general, shor1967generalized}. Метод зеркального спуска возник для безусловных задач в \cite{nemirovskii1979efficient,nemirovsky1983problem} как аналог стандартного субградиентного метода с неевклидовым проектированием. Для условных задач аналог этого метода был предложен в \cite{nemirovsky1983problem} (см. также \cite{beck2010comirror}). Проблема адаптивного выбора шага без использования констант Липшица рассмотрена в \cite{bib_Nemirovski} для задач без ограничений, а также в \cite{beck2010comirror} для задач с функциональными ограничениями. Однако критерии остановки для указанных методов зеркального спуска оставались неадаптивными.

Недавно в \cite{bib_Adaptive} были предложены алгоритмы зеркального спуска как с адаптивным выбором шага, так и с адаптивным критерием остановки. При этом помимо случая липшицевых целевого функционала и функционального ограничения был предложен оптимальный метод для класса условных задач выпуклой минимизации со специальными условиями роста (условимся называть их нестандартными) целевого функционала. Например, в задачах с квадратичными функционалами мы сталкиваемся с ситуацией, когда такой функционал не удовлетворяет обычному свойству Липшица (или константа Липшица достаточно большая), но градиент удовлетворяет условию Липшица. Для задач такого типа в (\cite{bib_Adaptive}, п. 3.3) был предложен адаптивный алгоритм зеркального спуска на базе идеологии \cite{bib_Nesterov,bib_Nesterov2016}.

В предлагаемой статье предложены новые модификации адаптивных методов зеркального спуска из (\cite{bib_Adaptive}, п. 3.1 и 3.3) для задач выпуклой минимизации в случае, когда имеется несколько функциональных ограничений. Отдельно рассмотрены методы для случаев липшицевости целевого функционала (пункт \ref{SectLip}), а также в случае нестандартных условий роста целевого функционала (пункт \ref{SectLipGrad}). Предлагаемые модификации позволяет сэкономить время работы алгоритма за счёт рассмотрения не всех функциональных ограничений на непродуктивных шагах. Получены оценки на скорость сходимости рассматриваемых методов (теоремы \ref{th1} и \ref{th2}, следствия \ref{cor1} и \ref{cor2}). Оценки необходимого числа итераций \eqref{eq08} и \eqref{eqq08}, указывают на оптимальность предлагаемых методов с точки зрения нижних оракульных оценок \cite{nemirovsky1983problem}. В заключении работы приведены результаты численных экспериментов, иллюстрирующие преимущества предлагаемых нами методов для некоторых примеров.

\subsection{Постановка задачи и необходимые вспомогательные понятия}

Пусть $(E,||\cdot||)$~--- конечномерное нормированное векторное пространство и $E^*$~--- сопряженное пространство к $E$ со стандартной нормой:
$$||y||_*=\max\limits_x\{\langle y,x\rangle,||x||\leq1\},$$
где $\langle y,x\rangle$~--- значение линейного непрерывного функционала $y$ в точке $x \in E$.

Пусть $X\subset E$~--- замкнутое выпуклое множество. Рассмотрим набор выпуклых субдифференцируемых функционалов $f$ и $g_m:X\rightarrow\mathbb{R}$ для всякого $m = \overline{1, M}$. Также предположим, что все функционалы $g_m$ удовлетворяют условию Липшица с некоторой константой $M_g$:
\begin{equation}\label{eq1}
|g_m(x)-g_m(y)|\leqslant M_g||x-y||\;\forall x,y\in X, \; m = \overline{1, M}.
\end{equation}

Будем рассматривать следующий тип задач выпуклой оптимизации.
\begin{equation}\label{eq2}
 f(x) \rightarrow \min\limits_{x\in X},
\end{equation}
где
\begin{equation}
\label{problem_statement_g}
    g_m(x) \leqslant 0 \quad \forall m = \overline{1, M}.
\end{equation}

Сделаем предположение о разрешимости задачи \eqref{eq2}--\eqref{problem_statement_g}.

Для дальнейших рассуждений нам потребуются вспомогательные понятия (см., например, \cite{bib_Nemirovski}). Введём так называемую {\it прокс-функцию} $d : X \rightarrow \mathbb{R}$, обладающую свойством непрерывной дифференцируемости и $1$-сильной выпуклости относительно нормы $\lVert\cdot\rVert$, т.е.
$$\langle \nabla d(x) - \nabla d(y), x-y \rangle \geqslant \lVert x-y \rVert^2 \quad \forall x, y, \in X$$
и предположим, что $\min\limits_{x\in X} d(x) = d(0).$ Будем полагать, что существует такая константа $\Theta_0 > 0$, что
$d(x_{*}) \leqslant \Theta_0^2,$ где $x_*$~--- точное решение (\ref{eq2})--(\ref{problem_statement_g}).
Если имеется множество решений $X_*$, то мы предполагаем, что для константы $\Theta_0$
$$\min\limits_{x_* \in X_*} d(x_*) \leqslant \Theta_0^2.$$
Для всех $x, y\in X$ рассмотрим соответствующую дивергенцию Брегмана
$$V(x, y) = d(y) - d(x) - \langle \nabla d(x), y-x \rangle.$$

В зависимости от постановки конкретной задачи возможны различные подходы к определению прокс-структуры задачи и соответствующей дивергенции Брегмана: евклидова, энтропийная и многие другие (см., например, \cite{bib_Nemirovski}). Стандартно определим оператор проектирования
$$\mathrm{Mirr}_x (p) = \arg\min\limits_{u\in X} \big\{ \langle p, u \rangle + V(x, u) \big\} \; \text{ для всяких  }x\in X \text{ и }p\in E^*.$$
Сделаем предположение о том, что оператор $\mathrm{Mirr}_x (p)$ легко вычислим.

Ясно, что вместо набора выпуклых функциональных ограничений $\{g_m(\cdot)\}_{m=1}^{M}$ можно рассмотреть одно ограничение $g: X \rightarrow \mathbb{R}$:
$$
g(x) = \max\limits_{m = \overline{1, M}} g_m(x), \quad |g(x)-g(y)|\leqslant M_g||x-y||\;\forall \quad x,y\in X.
$$

Для условных задач с целевым функционалом $f$ одним выпуклым субдифференцируемым функциональным ограничением в работе \cite{bib_Adaptive} предложено два метода зеркального спуска. Сходимость первого из них доказана для случая липшицевости целевого функционала (см. \cite{bib_Adaptive}, п. 3.1), сходимость же второго обоснована, в частности, в предположении липшицевости градиента $f$ (см. \cite{bib_Adaptive}, п. 3.3). Напомним эти методы.

\begin{algorithm}
\caption{Адаптивный зеркальный спуск (липшицев целевой функционал)}
\label{alg1}
\begin{algorithmic}[1]
\REQUIRE $\varepsilon>0,\Theta_0: \,d(x_*)\leqslant\Theta_0^2$
\STATE $x^0=argmin_{x\in X}\,d(x)$
\STATE $I=:\emptyset$
\STATE $N\leftarrow0$
\REPEAT
    \IF{$g(x^N)\leqslant\varepsilon$}
        \STATE $M_N=||\nabla f(x^N)||_*$, $h_N=\frac{\varepsilon}{M_N^2}$
        \STATE $x^{N+1}=Mirr_{x^N}(h_N\nabla f(x^N))\;\text{// \emph{"продуктивные шаги"}}$
        \STATE $N\rightarrow I$
    \ELSE
        \STATE $M_N=||\nabla g(x^N)||_*$, $h_N=\frac{\varepsilon}{M_N^2}$
        \STATE $x^{N+1}=Mirr_{x^N}(h_N\nabla g(x^N))\;\text{// \emph{"непродуктивные шаги"}}$
    \ENDIF
    \STATE $N\leftarrow N+1$
\UNTIL{$\sum\limits_{j=0}^{N-1}\frac{1}{M_j^2}\geqslant2\frac{\Theta_0^2}{\varepsilon^2}$}
\ENSURE $\bar{x}^N:=\frac{\sum\limits_{k\in I}x^kh_k}{\sum\limits_{k\in I}h_k}$
\end{algorithmic}
\end{algorithm}

\begin{algorithm}
\caption{Адаптивный зеркальный спуск (нестандартные условия роста)}
\label{alg2}
\begin{algorithmic}[1]
\REQUIRE $\varepsilon>0,\Theta_0: \,d(x_*)\leqslant\Theta_0^2$
\STATE $x^0=argmin_{x\in X}\,d(x)$
\STATE $I=:\emptyset$
\STATE $N\leftarrow0$
\REPEAT
    \IF{$g(x^N)\leqslant\varepsilon$}
        \STATE $h_N\leftarrow\frac{\varepsilon}{||\nabla f(x^N)||_{*}}$
        \STATE $x^{N+1}\leftarrow Mirr_{x^N}(h_N\nabla f(x^N))\;\text{// \emph{"продуктивные шаги"}}$
        \STATE $N\rightarrow I$
    \ELSE
        \STATE // \emph{$(g(x^N)>\varepsilon)$}
        \STATE $h_N\leftarrow\frac{\varepsilon}{||\nabla g(x^N)||_{*}^2}$
        \STATE $x^{N+1}\leftarrow Mirr_{x^N}(h_N\nabla g(x^N))\;\text{// \emph{"непродуктивные шаги"}}$
    \ENDIF
    \STATE $N\leftarrow N+1$
\UNTIL{$\Theta_0^2 \leqslant \frac{\varepsilon^2}{2}\left(|I|+\sum\limits_{k\not\in I}\frac{1}{||\nabla g(x^k)||_{*}^2}\right)$}
\ENSURE $\bar{x}^N:=argmin_{x^k,\;k\in I}\,f(x^k)$
\end{algorithmic}
\end{algorithm}

\begin{remark}
В обоих методах в ходе работы необходимо делить на нормы субградиентов целевого функционала $||\nabla f(x^k)||_{*}^2$ или ограничения $||\nabla g(x^k)||_{*}^2$. В связи с этим прокомментируем ситуацию, когда $\nabla f(x^k) = 0$ или $\nabla g(x^k) = 0$. Если верно $\nabla f(x^k) = 0$, то ясно, что $x^k$~--- точное решение задачи минимизации $f$ вне зависимости от ограничений и в этом случае работу алгоритма нужно останавливать. Условимся не оговаривать это обстоятельство отдельно в листингах алгоритмов. Если же $\nabla g(x^k) = 0$, то $x^k$~--- точка глобального минимума $g$. В сочетании с условием $g(x^k) > \varepsilon$ это означает, что $g(x) > \varepsilon > 0$ на всей области определения и тогда поставленная задача просто неразрешима. Ввиду изложенных обстоятельств условимся здесь и всюду далее полагать, что $\nabla f(x^k) \neq 0$ и $\nabla g(x^k) \neq 0$.
\end{remark}

В настоящей работе мы покажем, что на <<непродуктивных>> шагах ($k \not\in I$) можно вместо субградиента ограничения max-типа $g(x) = \max\limits_{m = \overline{1, M}} g_m(x)$ использовать субградиент любого из функционалов $g_m$, для которого верно $g_m(x^k) > \varepsilon$.

\subsection{Модификация алгоритма зеркального спуска в случае липшицевости целевого функционала}\label{SectLip}

Будем предполагать, что целевой функционал $f:X\rightarrow\mathbb{R}$ удовлетворяет условию Липшица:
\begin{equation}\label{eqf1}
|f(x)-f(y)|\leqslant M_f||x-y||\quad \forall x,y\in X.
\end{equation}
При таком соглашении мы рассмотрим алгоритм \ref{alg3} для задачи (\ref{eq2})--(\ref{problem_statement_g}).

\begin{algorithm}
\caption{Модифицированный адаптивный зеркальный спуск (липшицев целевой функционал)}
\label{alg3}
\begin{algorithmic}[1]
\REQUIRE $\varepsilon>0,\Theta_0: \,d(x_*)\leqslant\Theta_0^2$
\STATE $x^0=argmin_{x\in X}\,d(x)$
\STATE $I=:\emptyset$
\STATE $N\leftarrow0$
\REPEAT
    \IF{$g(x^N)\leqslant\varepsilon$}
        \STATE $M_N=||\nabla f(x^N)||_*$
        \STATE $h_N=\frac{\varepsilon}{M_N^2}$, $x^{N+1}=Mirr_{x^N}(h_N\nabla f(x^N))\;\text{// \emph{"продуктивные шаги"}}$
        \STATE $N\rightarrow I$
    \ELSE
        \STATE // \emph{$(g_{m(N)}(x^N)>\varepsilon)\;\text{для некоторого}\;  m(N)\in \{1,\ldots,M\}$}
        \STATE $M_N=||\nabla g_{m(N)}(x^N)||_*$, $h_N=\frac{\varepsilon}{M_N^2}$
        \STATE $x^{N+1}=Mirr_{x^N}(h_N\nabla g_{m(N)}(x^N))\;\text{// \emph{"непродуктивные шаги"}}$
    \ENDIF
    \STATE $N\leftarrow N+1$
\UNTIL{$\sum\limits_{j=0}^{N-1}\frac{1}{M_j^2}\geqslant2\frac{\Theta_0^2}{\varepsilon^2}$}
\ENSURE $\bar{x}^N:=\frac{\sum\limits_{k\in I}x^kh_k}{\sum\limits_{k\in I}h_k}$
\end{algorithmic}
\end{algorithm}

Напомним одно известное утверждение (см., например \cite{bib_Nemirovski}).
\begin{lemma}\label{lem1}
Пусть $f:X\rightarrow\mathbb{R}$~--- выпуклый субдифференцируемый функционал. Для произвольного $y \in X$ и некоторого $h > 0$ положим:
$$
z=Mirr_{y}(h\nabla f(y)).
$$
Тогда для произвольного $x\in X$
\begin{equation}\label{eq7}
h\langle\nabla f(y), y-x\rangle\leqslant\frac{h^2}{2}||\nabla f(y)||_*^2 + V(y,x) - V(z,x).
\end{equation}
\end{lemma}

\begin{theorem}\label{th1}
Пусть $\varepsilon > 0$~--- фиксированное число и выполнен критерий остановки алгоритма \ref{alg3}. Тогда $ \bar{x}^{N}$ есть $\varepsilon$-решение задачи \eqref{eq2}-- \eqref{problem_statement_g}:
$$
f(\bar{x}^N) - f^* \leqslant \varepsilon, \quad g(\bar{x}^N) \leqslant \varepsilon.
$$
При этом алгоритм \ref{alg3} работает не более
\begin{equation}\label{eq08}
N=\left\lceil\frac{2\max\{M_f^2,M_g^2\}\Theta_0^2}{\varepsilon^2}\right\rceil
\end{equation}
итераций.
\end{theorem}

\proof Мы отправляемся от (\cite{bib_Adaptive}, п. 3.1). Согласно определению $\bar{x}^N$ (см. алгоритм \ref{alg3}) и ввиду выпуклости целевого функционала $f$
\begin{equation}\label{th1_eq1}
        \sum\limits_{k\in I} h_k f(\bar{x}^N) - f(x_{*}) \leqslant \sum\limits_{k\in I} h_k \big( f(x^k) - f(x_{*}) \big).
\end{equation}

Далее, по лемме \ref{lem1}
$$
h_k \big(f(x^k) - f(x) \big) \leqslant \frac{h_k^2}{2} \lVert \nabla f (x^k) \rVert^2_* + V(x^k, x) - V(x^{k+1}, x) \quad \forall k \in I,
$$
$$
h_k \big(g_{m(k)}(x^k) - g_{m(k)}(x) \big) \leqslant \frac{h_k^2}{2} \lVert \nabla g_{m(k)}(x^k) \rVert^2_* + V(x^k, x) - V(x^{k+1}, x) \quad \forall k \not\in I
$$
и ввиду выбора величины шагов $h_k$ в алгоритме \ref{alg3} имеем:
\begin{gather*}
        \sum\limits_{k\in I} h_k \big( f(x^k) - f(x_{*}) \big) + \sum\limits_{k\not\in I} h_k \big(g_{m(k)}(x^k) - g_{m(k)}(x_{*}) \big) \leqslant \\
        \leqslant\sum\limits_{k\in I} \frac{h_k^2 ||\nabla f(x^k)||_*^2}{2} + \sum\limits_{k\not\in I} \frac{h_k^2 ||\nabla g_{m(k)}(x^k)||_*^2}{2} + \\
        +\sum\limits_{k = 0}^{N-1} \big( V(x^k, x_{*}) - V(x^{k+1}, x_{*}) \big) \leqslant \frac{\varepsilon}{2} \sum\limits_{k = 0}^{N-1} h_k + \Theta_0^2,
 \end{gather*}
где на непродуктивном шаге под $g_{m(k)}$ мы понимаем любое из ограничений, для которого верно неравенство $g_{m(k)}(x^k)>\varepsilon$. Поскольку для всякого $k \not\in I$
    $$g_{m(k)}(x^k) - g_{m(k)}(x_{*}) \geqslant g_{m(k)}(x^k) > \varepsilon,$$
то с учётом (\ref{th1_eq1}), мы имеем
    \begin{gather*}
        \sum\limits_{k\in I} h_k \big( f(\bar{x}^N) - f(x_{*}) \big) < \frac{\varepsilon}{2} \sum\limits_{k = 0}^{N-1} h_k - \varepsilon\sum\limits_{k\not\in I} h_k + \Theta_0^2 = \\
        =\varepsilon\sum\limits_{k\in I} h_k - \frac{\varepsilon^2}{2} \sum\limits_{k\in I} \frac{1}{||\nabla f(x^k)||_*^2} - \frac{\varepsilon^2}{2} \sum\limits_{k\not\in I} \frac{1}{||\nabla g_{m(k)}(x^k)||_*^2} + \Theta_0^2 \leqslant \varepsilon\sum\limits_{k\in I} h_k.
    \end{gather*}

Ввиду строгости последнего неравенства мы имеем, что множество продуктивных шагов $I$ непусто. Ясно, что для всякого $k\in I$ справедливо неравенство $g(x^k) \leqslant \varepsilon$. Тогда ввиду выпуклости $g$ верно
    \begin{gather*}
        \sum\limits_{k\in I} h_k g(\bar{x}^N) \leqslant \sum\limits_{k\in I} h_k g(x^k) \leqslant \varepsilon\sum\limits_{k\in I} h_k,
    \end{gather*}

Отметим также, что ввиду условий Липшица для целевого функционала и функционального ограничения (см. \eqref{eq1} и \eqref{eqf1}) на любой итерации работы алгоритма \ref{alg3} справедливы неравенства $||\nabla f(x^k)||_* \leqslant  M_f$ и $||\nabla g_{m(k)}(x^k)||_* \leqslant  M_g$. Поэтому критерий остановки алгоритма \ref{alg3} будет заведомо выполнен не более, чем после \eqref{eq08} итераций работы, что завершает доказательство теоремы.

\begin{remark}
Ввиду липшицевости и, вообще говоря, негладкости целевого функционала и ограничений оценка \eqref{eq08} на число итераций означает, что предложенный метод оптимален с точки зрения оракульных оценок \cite{bib_Nemirovski}: для достижения требуемой точности $\varepsilon$ решения задачи \eqref{eq2}--\eqref{problem_statement_g} достаточно $O\left(\frac{1}{\varepsilon^2}\right)$ итераций алгоритма \ref{alg3}.
\end{remark}

\subsection{Модифицированный алгоритм зеркального спуска в случае нестандартных условий роста целевого функционала}\label{SectLipGrad}

Теперь рассмотрим метод для класса условных задач выпуклой минимизации с нестандартными условиями роста целевого функционала. Например, в задачах с квадратичными целевыми функционалами мы сталкиваемся с ситуацией, когда такой функционал не удовлетворяет обычному свойству Липшица (или константа Липшица достаточно большая), но градиент удовлетворяет условию Липшица. Если целевой функционал $f$ не удовлетворяет свойству Липшица, то критерии остановки алгоритмов \ref{alg1} и \ref{alg3} не позволяют получить оценку вида \eqref{eq08} и обосновать оптимальность метода с точки зрения нижних оракульных оценок. Такие проблемы возникают и для более широкого класса уже негладких целевых функционалов
\begin{equation}\label{equiv_nonstand1}
f(x)=\max\limits_{1\leqslant i\leqslant m} f_i(x),
\end{equation}
где
\begin{equation}\label{equiv_nonstand2}
f_i(x)=\frac{1}{2}\langle A_ix,x\rangle-\langle b_i,x\rangle+\alpha_i,\;i=1,\ldots,m,
\end{equation}
в случае, когда $A_i$ ($i=1,\ldots,m$)~--- положительно определённые матрицы: $x^TA_ix\geqslant 0\ \forall x \in X$.

Отметим, что функционалы вида \eqref{equiv_nonstand1}--\eqref{equiv_nonstand2} возникают в задачах проектирования механических конструкций Truss Topology Design со взвешенными балками (см., например, презентацию \cite{bib_Nesterov2016}).

Для задач \eqref{eq2}--\eqref{problem_statement_g} с целевыми функционалами вида \eqref{equiv_nonstand1}--\eqref{equiv_nonstand2}  в (\cite{bib_Adaptive}, п. 3.3) был предложен адаптивный алгоритм типа зеркального спуска на базе идеологии \cite{bib_Nesterov,bib_Nesterov2016}. Мы рассмотрим модификацию этого метода в случае наличия нескольких функциональных ограничений. Аналогично (\cite{bib_Nesterov}, п. 3.2.2), для всякого субградиента $\nabla f(x)$ целевого функционала $f$ в точке $y \in X$, введём следующую вспомогательную величину
\begin{equation}
v_f(x, y)=\left\{
\begin{aligned}
&\left\langle\frac{\nabla f(x)}{\|\nabla f(x)\|_{*}},x-y\right\rangle, \quad &\nabla f(x) \ne 0\\
&0 &\nabla f(x) = 0\\
\end{aligned}
\right.,\quad x \in X.
\label{eq:vfDef}
\end{equation}

Рассмотрим следующий алгоритм адаптивного зеркального спуска для задач (\ref{eq2})--(\ref{problem_statement_g}).

\begin{algorithm}
\caption{Модифицированный адаптивный зеркальный спуск (нестандартные условия роста)}
\label{alg4}
\begin{algorithmic}[1]
\REQUIRE $\varepsilon>0,\Theta_0:\,d(x_*)\leqslant\Theta_0^2$
\STATE $x^0=argmin_{x\in X}\,d(x)$
\STATE $I=:\emptyset$
\STATE $N\leftarrow0$
\REPEAT
    \IF{$g(x^N)\leqslant\varepsilon$}
        \STATE $h_N\leftarrow\frac{\varepsilon}{||\nabla f(x^N)||_{*}}$
        \STATE $x^{N+1}\leftarrow Mirr_{x^N}(h_N\nabla f(x^N))\;\text{// \emph{"продуктивные шаги"}}$
        \STATE $N\rightarrow I$
    \ELSE
        \STATE // \emph{$(g_{m(N)}(x^N)>\varepsilon)\;\text{для некоторого}\; m(N)\in \{1,\ldots,M\}$}
        \STATE $h_N\leftarrow\frac{\varepsilon}{||\nabla g_{m(N)}(x^N)||_{*}^2}$
        \STATE $x^{N+1}\leftarrow Mirr_{x^N}(h_N\nabla g_{m(N)}(x^N))\;\text{// \emph{"непродуктивные шаги"}}$
    \ENDIF
    \STATE $N\leftarrow N+1$
\UNTIL{$\Theta_0^2 \leqslant \frac{\varepsilon^2}{2}\left(|I|+\sum\limits_{k\not\in I}\frac{1}{||\nabla g_{m(k)}(x^k)||_{*}^2}\right)$}
\ENSURE $\bar{x}^N:=argmin_{x^k,\;k\in I}\,f(x^k)$
\end{algorithmic}
\end{algorithm}

Справедлива следующая
\begin{theorem}\label{th2}
Пусть $\varepsilon > 0$~--- фиксированное число и выполнен критерий остановки алгоритма \ref{alg4}. Тогда
\begin{equation}\label{eq09}
\min\limits_{k \in I} v_f(x^k,x_*)<\varepsilon.
\end{equation}
Отметим, что алгоритм \ref{alg4} работает не более
\begin{equation}\label{eqq08}
N=\left\lceil\frac{2\max\{1, M_g^2\}\Theta_0^2}{\varepsilon^2}\right\rceil
\end{equation}
итераций.
\end{theorem}
\proof Мы отправляемся от (\cite{bib_Adaptive}, п. 3.3). Пусть $[N]=\{k\in\overline{0,N-1}\}$, $J=[N]\setminus I$~--- набор номеров непродуктивных шагов.

1) Для продуктивных шагов ввиду леммы \ref{lem1} (см. \eqref{eq7}) имеем, что $$h_k\langle\nabla f(x^k), x^k-x\rangle\leqslant\frac{h_k^2}{2}||\nabla f(x^k)||_*^2+V(x^k,x)-V(x^{k+1},x).$$

Примем во внимание $\frac{h_k^2}{2}||\nabla f(x^k)||_*^2=\frac{\varepsilon^2}{2}$, имеем
\begin{equation}\label{eq001}
\begin{split}
h_k\langle\nabla f(x^k),x^k-x\rangle=\varepsilon\left\langle\frac{\nabla f(x^k)}{||\nabla f(x^k)||_*},\, x^k-x\right\rangle=\varepsilon v_f(x^k,x) \leqslant\\
\leqslant \frac{\varepsilon^2}{2} + V(x^k,x_*)-V(x^{k+1},x_*).
\end{split}
\end{equation}

2) Аналогично для «непродуктивных» шагов $k\in J$ (под $g_m(k)$ мы понимаем любое ограничение, для которого $g_m(k)>\varepsilon$) по лемме \ref{lem1}:
$$
h_k\left(g_{m(k)}(x^k)-g_{m(k)}(x_*)\right)\leqslant\frac{h_k^2}{2}||\nabla g_{m(k)}(x^k)||_*^2+V(x^k,x_*)-V(x^{k+1},x_*)=
$$
$$
=\frac{\varepsilon^2}{2||\nabla g_{m(k)}(x^k)||_*^2}+V(x^k,x_*)-V(x^{k+1},x_*).
$$

3) Из (\ref{eq1}) и (\ref{eq2}) при $x=x_*$ имеем:
\begin{equation}\label{eq10}
\varepsilon\cdot \sum\limits_{k\in I}v_f(x^k,x_*)+\sum\limits_{k\in J}\frac{\varepsilon^2 (g_{m(k)}(x^k)-g_{m(k)}(x_*))}{2||\nabla g_{m(k)}(x^k)||_*^2}\leqslant \frac{\varepsilon^2}{2} \cdot |I| +
\sum\limits_{k=0}^{N-1}(V(x^k,x_*)-V(x^{k+1},x_*)).
\end{equation}

Заметим, что для любого $k \in J$
\begin{equation}\label{eqiavgmk}
g_{m(k)}(x^k)-g_{m(k)}(x_*)\geqslant g_{m(k)}(x^k)>\varepsilon
\end{equation}
и ввиду
$$\sum\limits_{k=0}^{N-1}(V(x^k,x_*)-V(x^{k+1},x_*))\leqslant\Theta_0^2$$
неравенство (\ref{eq10}) может быть преобразовано следующим образом:
$$\varepsilon\sum\limits_{k\in I}v_f(x^k,x_*)\leqslant |I|\cdot\frac{\varepsilon^2}{2}+\Theta_0^2-\sum\limits_{k \in J}\frac{\varepsilon^2}{2||\nabla g_{m(k)}(x^k)||_*^2},$$
$$\sum\limits_{k\in I}v_f(x^k,x_*)\geqslant|I|\min\limits_{k\in I}v_f(x^k,x_*).$$

Таким образом,
$$\varepsilon \cdot \min\limits_{k \in I} v_f(x^k,x_*)\cdot |I| \leqslant \frac{\varepsilon^2}{2}\cdot|I|+\Theta_0^2-\sum\limits_{k \in J}\frac{\varepsilon^2}{2||\nabla g_{m(k)}(x^k)||_*^2}\leqslant \varepsilon^2\cdot|I|,$$
откуда
\begin{equation}\label{eq12}
\min\limits_{k \in I} v_f(x^k,x_*) \leqslant \varepsilon.
\end{equation}

В завершении покажем, что $|I|\neq 0$. Предположим обратное: $|I|=0\Rightarrow|J|=N$, т.е. все шаги непродуктивны. Тогда с учётом \eqref{eqiavgmk} получаем, что
$$
h_k(g_{m(k)}(x^k)-g_{m(k)}(x_*)) > \frac{\varepsilon^2}{\|\nabla g_{m(k)}(x^k)\|_*^2}
$$
и
$$
\sum\limits_{k=0}^{N-1} h_k(g_{m(k)}(x^k)-g_{m(k)}(x_*))\leqslant\sum\limits_{k=0}^{N-1}\frac{\varepsilon^2}{2\|\nabla g_{m(k)}(x^k)\|_*^2}+\Theta_0^2\leqslant
\sum\limits_{k=0}^{N-1}\frac{\varepsilon^2}{\|\nabla g_{m(k)}(x^k)\|_*^2}.
$$

Итак, получили противоречие. Это означает, что $|I|\neq 0$.

Отметим, что ввиду условия Липшица для функционального ограничения на любой итерации работы алгоритма \ref{alg4} справедливо неравенство $||\nabla g_{m(k)}(x^k)||_* \leqslant  M_g$. Поэтому критерий остановки алгоритма \ref{alg4} будет заведомо выполнен не более, чем после \eqref{eqq08} итераций работы. Теорема доказана.

Теперь покажем, как можно оценить скорость сходимости предлагаемого метода. Для этого полезно следующее вспомогательное утверждение (\cite{bib_Nesterov}, лемма 3.2.1). Напомним, что $x_*$~--- решение задачи \eqref{eq2}--\eqref{problem_statement_g}.
\begin{lemma}
Введем следующую функцию:
\begin{equation}\label{eq13}
\omega(\tau)=\max\limits_{x\in X}\{f(x)-f(x_*):||x-x_*||\leqslant\tau\},
\end{equation} где $\tau$ - положительное число.
Тогда для всякого $y \in X$
\begin{equation}\label{eq_lemma}
f(y) - f(x_*) \leqslant \omega(v_f(y,x_*)).
\end{equation}
\end{lemma}

На базе предыдущего утверждения и теоремы \ref{th2} можно оценить скорость сходимости алгоритма \ref{alg4} для дифференцируемого целевого функционала $f$ с градиентом, удовлетворяющим условию Липшица:
\begin{equation}\label{eqlipgrad}
||\nabla f(x)-\nabla f(y)||_*\leqslant L||x-y|| \quad \forall x,y\in X.
\end{equation}

Используя следующее известное неравенство (см., например \cite{bib_Nesterov})
$$f(x)\leqslant f(x_*)+||\nabla f(x_*)||_*||x-x_*||+\frac{1}{2}L||x-x_*||^2,$$
мы можем получить, что
$$\min\limits_{k\in I}f(x^k)-f(x_*)\leqslant\min\limits_{k\in I} \left\{||\nabla f(x_*)||_*||x^k-x_*||+\frac{1}{2}L||x^k-x_*||^2\right\}.$$
Далее, справедливы оценки:
$$
f(x)-f(x_*)\leqslant \varepsilon \cdot ||\nabla f(x_*)||_* + \frac{1}{2} L\varepsilon^2.
$$

Поэтому справедливо

\begin{corollary}\label{cor1}
Пусть $f$ дифференцируема на $X$ и верно \eqref{eqlipgrad}. Тогда после остановки алгоритма \ref{alg4} верна оценка: $$\min\limits_{1\leqslant k\leqslant N}f(x^k)-f(x_*)\leqslant\varepsilon_f+\frac{L\varepsilon^2}{2},
$$ где $$\varepsilon_f=\varepsilon \cdot ||\nabla f(x_*)||_*.$$
\end{corollary}

Также можно рассмотреть специальный класс негладких целевых функционалов.

\begin{corollary}\label{cor2}
Пусть $f(x) = \max\limits_{i = \overline{1, m}} f_i(x)$, где $f_i$ дифференцируема для всякого $x \in X$ и
$$
||\nabla f_i(x)-\nabla f_i(y)||_*\leqslant L_i||x-y|| \quad \forall x,y\in X.
$$
Тогда после остановки алгоритма \ref{alg4} верна оценка: $$\min\limits_{0\leqslant k\leqslant N}f(x^k)-f(x_*)\leqslant\varepsilon_f+\frac{L\varepsilon^2}{2},
$$ где $$\varepsilon_f=\varepsilon \cdot ||\nabla f(x_*)||_*, \quad L = \max\limits_{i = \overline{1, m}} L_i.$$
\end{corollary}

\begin{remark}
Ввиду липшицевости и, вообще говоря, негладкости функциональных ограничений оценка \eqref{eqq08} на число итераций означает, что предложенный метод оптимален с точки зрения оракульных оценок \cite{bib_Nemirovski}: для достижения требуемой точности $\varepsilon$ решения задачи \eqref{eq2}--\eqref{problem_statement_g} для рассмотренного в данном разделе статьи класса целевых функционалов достаточно $O\left(\frac{1}{\varepsilon^2}\right)$ итераций алгоритма \ref{alg4}. Отметим также, что к рассмотренным классам задач \eqref{eq2}--\eqref{problem_statement_g} с целевыми функционалами вида \eqref{equiv_nonstand1}--\eqref{equiv_nonstand2} применимы и рассмотренные ранее алгоритмы \ref{alg1} и \ref{alg3}. Однако невыполнение, вообще говоря, условия Липшица для $f$ не позволяет обосновать оптимальность алгоритмов \ref{alg1} и \ref{alg3} для целевых функционалов вида \eqref{equiv_nonstand1}--\eqref{equiv_nonstand2}. Точнее говоря, возможны ситуации, когда на продуктивных шагах нормы (суб)градиентов целевого функционала $\|\nabla f(x^k)\|_*$ будут достаточно большими и это будет мешать быстрому достижению критерия остановки алгоритмов \ref{alg1} и \ref{alg3}. В таком случае алгоритмов \ref{alg2} и \ref{alg4} могут работать быстрее, что показано ниже на некоторых конкретных примерах в численных экспериментах.
\end{remark}

\subsection{Численные эксперименты}

Для демонстрации преимуществ алгоритмов \ref{alg3} и \ref{alg4} по сравнению с алгоритмами \ref{alg1} и \ref{alg2} соответственно был проведен ряд тестов. Рассматриваемые целевые функционалы $f(x)$ и ограничения $g_m(x)$ ($m = 1, 2, 3, \ldots, 10$) указаны в таблицах 1 (более простые целевые функционалы) и 3 (негладкие целевые функционалы max-типа). Отметим, что мы рассматриваем функционалы 10 переменных, стандартную евклидову прокс-структуру, начальную точку $x^0 = (1, 1, 1, \ldots, 1)$ и $\Theta_0=3$ (можно проверить, что для всех примеров одной из оптимальных точек будет $x_* = (0, 0, 0, \ldots, 0)$), точность $\varepsilon = 0.05$.

Результаты выполнения алгоритмов представлены в сравнительных таблицах 2 и 4. Приводится количество итераций и время (указано в секундах) работы алгоритмов \ref{alg1} и \ref{alg2}, а также соответствующих модификаций. Все вычисления были произведены с помощью CPython 3.6.4 на компьютере с 3-ядерным процессором AMD Athlon II X3 450 с тактовой частотой 803,5 МГц на каждое ядро. ОЗУ компьютера составляла 8 Гб.

\begin{table}[]
\centering
\caption{"Входные данные"}
\label{tab1}
\begin{tabular}{|c|c|l|l|l|}
\hline
\multirow{2}{*}{} & \multicolumn{4}{c|}{Алгоритмы 1--4}                                                                                                                                                                                                                                                                                                                                                                                                                                                                                                                                                                                                                                                                                                                                                                                                \\ \cline{2-5}
                  & \multicolumn{4}{c|}{$f(x)$}                                                                                                                                                                                                                                                                                                                                                                                                                                                                                                                                                                                                                                                                                                                                                                                                        \\ \hline
Пр. 1             & \multicolumn{4}{c|}{$\sqrt{0.1\left(\sum\limits_{i=1}^{10}x_i^2+\sum\limits_{i=1}^9x_ix_{i+1}\right)}$}                                                                                                                                                                                                                                                                                                                                                                                                                                                                                                                                                                                                                                                                                                                            \\ \hline
Пр. 2             & \multicolumn{4}{c|}{$\sum\limits_{i=1}^{10}x_i^2-x_1x_2+x_3-x_8+x_9x_{10}$}                                                                                                                                                                                                                                                                                                                                                                                                                                                                                                                                                                                                                                                                                                                                                        \\ \hline
Пр. 3             & \multicolumn{4}{c|}{$\sum\limits_{i=1}^{10}5^ix_i^2$}                                                                                                                                                                                                                                                                                                                                                                                                                                                                                                                                                                                                                                                                                                                                                                              \\ \hline
\multirow{2}{*}{} & \multicolumn{4}{c|}{$g_m(x),\;m=\overline{1,10}$}                                                                                                                                                                                                                                                                                                                                                                                                                                                                                                                                                                                                                                                                                                                                                                                  \\ \cline{2-5}
                  & \multicolumn{4}{c|}{\begin{tabular}[c]{@{}c@{}}$x_1+20x_2+30x_3+40x_4+50x_5+60x_6+70x_7+80x_8+90x_9+100x_{10},$\\ $x_1+120x_2+130x_3+140x_4+150x_5+160x_6+170x_7+180x_8+190x_9+200x_{10},$\\ $x_1+220x_2+230x_3+240x_4+250x_5+260x_6+270x_7+280x_8+290x_9+300x_{10},$\\ $x_1+320x_2+330x_3+340x_4+350x_5+360x_6+370x_7+380x_8+390x_9+400x_{10},$\\ $x_1+420x_2+430x_3+440x_4+450x_5+460x_6+470x_7+480x_8+490x_9+500x_{10},$\\ $x_1+520x_2+530x_3+540x_4+550x_5+560x_6+570x_7+580x_8+590x_9+600x_{10},$\\ $x_1+620x_2+630x_3+640x_4+650x_5+660x_6+670x_7+680x_8+690x_9+700x_{10},$\\ $x_1+720x_2+730x_3+740x_4+750x_5+760x_6+770x_7+780x_8+790x_9+800x_{10},$\\ $x_1+820x_2+830x_3+840x_4+850x_5+860x_6+870x_7+880x_8+890x_9+900x_{10},$\\ $x_1+920x_2+930x_3+940x_4+950x_5+960x_6+970x_7+980x_8+990x_9+1\,000x_{10}$\end{tabular}} \\ \hline
\end{tabular}
\end{table}

\begin{table}[]
\centering
\caption{"Сравнение результатов работы алгоритмов"}
\label{tab2}
\begin{tabular}{|c|c|c|c|c|}
\hline
\multirow{2}{*}{} & Итерации          & Время       & Итерации          & Время       \\ \cline{2-5}
                  & \multicolumn{2}{c|}{Алгоритм 1} & \multicolumn{2}{c|}{Алгоритм 3} \\ \hline
Пр. 1             & 730\,829          & 133         & 261\,800          & 40          \\ \hline
Пр. 2             & 1\,638\,946       & 262         & 453\,580          & 30          \\ \hline
Пр. 3             & >$10^7$           & >500        & >$10^7$           & >500        \\ \hline
                  & \multicolumn{2}{c|}{Алгоритм 2} & \multicolumn{2}{c|}{Алгоритм 4} \\ \hline
Пр. 1             & --                & --          & --                & --          \\ \hline
Пр. 2             & 1\,584\,616       & 300         & 1\,434\,006       & 156         \\ \hline
Пр. 3             & 184\,706          & 124         & 89\,940           & 110         \\ \hline
\end{tabular}
\end{table}

\begin{table}[]
\centering
\caption{"Входные данные"}
\label{tab3}
\begin{tabular}{|c|c|l|l|l|}
\hline
\multirow{2}{*}{\textbf{}} & \multicolumn{4}{c|}{Алгоритмы 1--4}                                                                                                                                                                                                                                                                                                                                                                                                                                                                                                                                                                                                                                                                                                                                                                                                \\ \cline{2-5}
                           & \multicolumn{4}{c|}{$f(x)$}                                                                                                                                                                                                                                                                                                                                                                                                                                                                                                                                                                                                                                                                                                                                                                                                        \\ \hline
Пр. 4                      & \multicolumn{4}{c|}{\begin{tabular}[c]{@{}c@{}}$\max\{0.1|x_1+x_2+x_3|+1,\;0.01|x_4+2x_5+x_6|+2,$\\ $0.001|x_7+3x_8+4x_9+10x_{10}|+5\}$\end{tabular}}                                                                                                                                                                                                                                                                                                                                                                                                                                                                                                                                                                                                                                                                              \\ \hline
Пр. 5                      & \multicolumn{4}{c|}{\begin{tabular}[c]{@{}c@{}}$\max\{x_1^2,\;10x_2^2,\;50x_3^2,\;100x_4^2,\;200x_5^2,\;400x_6^2,$\\ $800x_7^2,\;1\,000x_8^2,\;5\,000x_9^2,\;10\,000x_{10}^2\}$\end{tabular}}                                                                                                                                                                                                                                                                                                                                                                                                                                                                                                                                                                                                                                      \\ \hline
Пр. 6                      & \multicolumn{4}{c|}{\begin{tabular}[c]{@{}c@{}}$\max\{x_1+2x_2+3x_3,\;x_3+4x_4+6x_5,\;x_4+3x_5+6x_6+7x_7,$\\ $5x_7+8x_8+9x_9,\;x_1+10x_{10}\}$\end{tabular}}                                                                                                                                                                                                                                                                                                                                                                                                                                                                                                                                                                                                                                                                       \\ \hline
\multirow{2}{*}{\textbf{}} & \multicolumn{4}{c|}{$g_m(x),\;m=\overline{1,10}$}                                                                                                                                                                                                                                                                                                                                                                                                                                                                                                                                                                                                                                                                                                                                                                                  \\ \cline{2-5}
                           & \multicolumn{4}{c|}{\begin{tabular}[c]{@{}c@{}}$x_1+20x_2+30x_3+40x_4+50x_5+60x_6+70x_7+80x_8+90x_9+100x_{10},$\\ $x_1+120x_2+130x_3+140x_4+150x_5+160x_6+170x_7+180x_8+190x_9+200x_{10},$\\ $x_1+220x_2+230x_3+240x_4+250x_5+260x_6+270x_7+280x_8+290x_9+300x_{10},$\\ $x_1+320x_2+330x_3+340x_4+350x_5+360x_6+370x_7+380x_8+390x_9+400x_{10},$\\ $x_1+420x_2+430x_3+440x_4+450x_5+460x_6+470x_7+480x_8+490x_9+500x_{10},$\\ $x_1+520x_2+530x_3+540x_4+550x_5+560x_6+570x_7+580x_8+590x_9+600x_{10},$\\ $x_1+620x_2+630x_3+640x_4+650x_5+660x_6+670x_7+680x_8+690x_9+700x_{10},$\\ $x_1+720x_2+730x_3+740x_4+750x_5+760x_6+770x_7+780x_8+790x_9+800x_{10},$\\ $x_1+820x_2+830x_3+840x_4+850x_5+860x_6+870x_7+880x_8+890x_9+900x_{10},$\\ $x_1+920x_2+930x_3+940x_4+950x_5+960x_6+970x_7+980x_8+990x_9+1\,000x_{10}$\end{tabular}} \\ \hline
\end{tabular}
\end{table}

\begin{table}[]
\centering
\caption{"Сравнение результатов работы алгоритмов"}
\label{tab4}
\begin{tabular}{|c|c|c|c|c|}
\hline
\multirow{2}{*}{\textbf{}} & Итерации         & Время        & Итерации         & Время        \\ \cline{2-5}
                           & \multicolumn{2}{c|}{Алгоритм 1} & \multicolumn{2}{c|}{Алгоритм 3} \\ \hline
Пр. 4                      & 172\,821         & 24           & 17\,255          & 1            \\ \hline
Пр. 5                      & >$10^6$          & >500         & >$10^6$          & >500         \\ \hline
Пр. 6                      & >$10^6$          & >500         & >$10^6$          & >500         \\ \hline
                           & \multicolumn{2}{c|}{Алгоритм 2} & \multicolumn{2}{c|}{Алгоритм 4} \\ \hline
Пр. 4                      & --               & --           & --               & --           \\ \hline
Пр. 5                      & 182\,993         & 106          & 66\,095          & 79           \\ \hline
Пр. 6                      & 180\,020         & 101          & 24\,454          & 78           \\ \hline
\end{tabular}
\end{table}

Как видим, предлагаемые нами модификации (алгоритм \ref{alg3}~--- модификация алгоритма \ref{alg1} и алгоритм \ref{alg4}~--- модификация алгоритма \ref{alg2}) могут существенно сокращать как необходимое для достижения нужного качества решения количество необходимых итераций, так и время работы алгоритмов для задач с указанными ограничениями. Особенно наглядно это видно в примерах 1, 2 и 4 для алгоритмов \ref{alg1} и \ref{alg3}. Также заметные преимущества даёт алгоритм \ref{alg4} перед алгоритмом \ref{alg2} в примерах 3, 5 и 6.

Это связано с тем, что адаптивные критерии остановки алгоритмов \ref{alg1}--\ref{alg4}, как легко видеть, определяют зависимость остановки алгоритмов от величины норм (суб)градиентов ограничений на непродуктивных шагах. Чем меньше эти нормы, тем скорее будут выполнены критерии остановки. Точнее говоря, замена (суб)градиента max-ограничения $g(\cdot) = \max\limits_{m \in \overline{1, M}} g_m(\cdot)$ (если применять алгоритм \ref{alg1} или \ref{alg2}) на какой-либо (суб)градиент $g_m(\cdot)$ (если применять алгоритм \ref{alg3} или \ref{alg4}) может ускорить выполнение критерия остановки метода.

Показано также, что в случае целевых функционалов $f$ вида \eqref{equiv_nonstand1}--\eqref{equiv_nonstand2} ввиду того, что величины $\|\nabla f(x^k)\|_*$ могут быть достаточно большими, алгоритмы \ref{alg2} и \ref{alg4} могут работать существенно быстрее алгоритмов \ref{alg1} и \ref{alg3}. Это показывают примеры 3, 5 и 6. Однако, если нормы $\|\nabla f(x^k)\|_*$ не особо велики, то алгоритмы \ref{alg1} и \ref{alg3} могут в конкретных примерах работать быстрее алгоритмов \ref{alg2} и \ref{alg4} и для квадратичных функционалов (см. пример 2).

\subsection{Заключительные замечания}

В заключении отметим, что основные результаты работы (сходимость и оптимальность предложенных методов с точки зрения нижних оракульных оценок) не зависят от выбора на непродуктивном шаге ограничения $g_m(\cdot)$, на котором $g_m(x^k) > \varepsilon$.

Однако замена (суб)градиента max-ограничения $g(\cdot) = \max\limits_{m \in \overline{1, M}} g_m(\cdot)$ (если применять алгоритм \ref{alg1} или \ref{alg2}) на какой-либо (суб)градиент $g_m(\cdot)$ (если применять алгоритм \ref{alg3} или \ref{alg4}) может привести к увеличению нормы рассматриваемого на непродуктивном шаге (суб)градиента. Чтобы этого избежать, нужно на каждом непродуктивном шаге $k$ среди ограничений $g_m(\cdot)$, для которых $g_m(x^k) > \varepsilon$, находить ограничение с наименьшей нормой (суб)градиента $\|\nabla g_{m}(x^k)\|$. Ясно, что операция минимизации $\|\nabla g_{m}(x^k)\|$ по всем подходящим ограничениям на каждой итерации может привести к увеличению времени работы алгоритма и в результате будут не очевидны преимущества перед исходными алгоритмами \ref{alg1} и \ref{alg2}.

Однако отметим, что предлагаемые нами модификации алгоритмов \ref{alg1} и \ref{alg2} выгодны в случае, если возможно априорно независимо от выбора точки упорядочить по неубыванию нормы (суб)градиентов функциональных ограничений. Это легко сделать, например, если задача имеет аффинные функциональные ограничения. Ясно, что существуют и другие примеры подходящих задач с ограничениями. Если такая сортировка возможна, то в алгоритмах \ref{alg3} и \ref{alg4} на непродуктивных шагах достаточно использовать первое ограничение $g_m(\cdot)$, для которого $g_m(x^k) > \varepsilon$. Представляется, что дальнейшая разработка предложенных методов для условных задач с различной структурой~--- интересная задача на будущее.

\subsection*{Благодарности} Авторы выражают огромную признательность Гасникову Александру Владимировичу и Двуреченскому Павлу Евгеньевичу, а также неизвестному рецензенту за полезные обсуждения и пожелания.

\noindent Федор Сергеевич Стонякин\\
кандидат физ.-мат. наук, доцент\\
Крымский федеральный университет им. В.~И.~Вернадского\\

\noindent Mohammad S. Alkousa\\
аспирант\\
Московский физико-технический институт (государственный университет)\\

\noindent Алексей Николаевич Степанов \\
соисполнитель работ по гранту Президента РФ для молодых кандидатов наук, код МК-176.2017.1\\
Крымский федеральный университет им. В.~И.~Вернадского\\

\noindent Максим Александрович Баринов\\
студент\\
Крымский федеральный университет им. В.~И.~Вернадского\\

\end{document}